\documentclass[runningheads]{llncs}
\usepackage{graphicx}
\usepackage{amsfonts}
\usepackage{amsmath}
\usepackage{amssymb,latexsym}
\usepackage{color}
\usepackage{graphicx}
\usepackage{cite}
\usepackage{hyperref} 
\hypersetup{colorlinks=true, linkcolor=blue,  anchorcolor=blue,  
citecolor=blue, filecolor=blue, menucolor=blue, pagecolor=blue,  
urlcolor=blue, pdftitle={DPn}}

\begin{document}
\title{About the Study of the $n$-dimensional Boolean Cube in the Bachelor's Programs in Computer Science
%%\thanks{Partially supported by the Research Fund of the University of Veliko Tarnovo, Bulgaria.}
}
\titlerunning{About the Study of the Boolean Cube in the Undergraduate Programs in CS}
\author{Valentin Bakoev}
\authorrunning{V. Bakoev}
\institute{``St. Cyril and St. Methodius'' University of Veliko Tarnovo, Bulgaria.
\email{v.bakoev@ts.uni-vt.bg}
}
\maketitle             
\begin{abstract}
	Here we present our arguments for a more in-depth study of the Boolean cube, which is one of the most important discrete structures. The article contains a case study that analyses how the Boolean cube has been included and explained in more than 80 different sources. However, organized material on the Boolean cube is lacking. We examine and show why the topic of the Boolean cube deserves to be studied in the course on Discrete Structures -- a basic part of the Computer Science curriculum. The benefits of mastering such knowledge and programming skills are pointed out. A sample lecture on the $n$-dimensional Boolean cube (including selected exercises, their answers, hints or solutions) is developed and discussed. It introduces, generalizes and relates many concepts from different subjects in the area of Discrete Mathematics and outside it. So the lecturers can use the sample lecture when teaching those subjects. In this way, all these concepts become more understandable, applicable and useful to the students. The lecture can be taught in a form of presentation or given to the students for ``Online self-learning", in order to write homework, or develop a course project related to the topic, etc. All supporting resources to the lecture are available for free use at \cite{VB_Lect}. 
\end{abstract}


\vspace{4mm}

An abridged version of this article has been published -- see Bakoev V., On the Study of the n-dimensional Boolean Cube in the Undergraduate Programs in Computer Science, 2021 International Conference Automatics and Informatics (ICAI), 2021, pp. 296-299, doi: 10.1109/ICAI52893.2021.9639496, \url{https://ieeexplore.ieee.org/abstract/document/9639496}.
\section{Introduction}
\label{Intro}
	The set $\{0,1\}^n$ of all $n$-dimensional binary vectors is known as \textit{$n$-dimensional Boolean cube} (\textit{hypercube}). It has many useful properties and applications in Computer science and in the practice of the programmers. It is well-known that the Boolean cube: 
\begin{itemize}
	\item Is an important concept in many subjects in the area of Discrete Mathematics -- for example, Combinatorics, Boolean Functions (hence in Computer Architectures), Coding Theory, Cryptology, Data Compression, Graph Theory, Theory of Algorithms, etc. Being their common part, the Boolean cube relates them.
	\item It and its properties are widely used in many areas outside Discrete Mathematics, such as Programming languages, Algorithms and Data Structures (especially in Combinatorial Algorithms, Cryptographic Algorithms, etc). A lot of techniques are based on these properties, which allow to create, analyze, illustrate or speed-up some algorithms.
\end{itemize}
 
		For more than 17 years, the Discrete Structures/Discrete Mathematics (DS/ DM) modules taught by the author have included a lecture on the $n$-dimensional Boolean cube (these modules are taught in 4 undergraduate programs). The lecture is placed after those for combinatorial configurations and precedes the study of Boolean functions. Our experience shows that studying this topic is very useful to the students, to the author (in teaching the remaining part of the DS/DM module, and also in teaching Algorithms and Data Structures, Cryptography, Competitive Programming), as well as to some of our colleagues in teaching other subjects. During the first several (of those 17) years we have been using a few selected textbooks. In the past several years we have examined how the $n$-dimensional Boolean cube is represented in numerous textbooks, books, lecture notes, etc., available to us. We could not find a whole topic (a part or a chapter) dedicated to it in these sources. These reasons have motivated us to discuss our ideas with colleagues and students, to represent them at a conference, to include such a chapter in the textbook \cite{VB2014}, and finally, to write this paper. By it, \textbf{we aim} to share our experience and to show why and how this topic can be included in the DS/DM module. Furthermore, if some lecturers decide to teach this topic, we can help them with the full text of a sample lecture that includes some selected problems and their hints, answers or solutions. All necessary files in .pdf and .tex formats are uploaded for free access at \cite{VB_Lect}. So, the teachers can easily modify the lecture and prepare a presentation in accordance with the content of their module and their educational objectives. If there is not enough time for including such a lecture in the DM/DS module, these resources could be pointed out to the students for self-study. They could support them when solving problems or working on projects part of their course. Due to the restrictions imposed by COVID'19, e-learning has become almost mandatory. The need for (free) educational resources is constantly growing. We hope that this article and the supporting electronic resources can help in teaching and studying the topic of the Boolean cube.
	
	We note that the lecture can be represented from an algebraic point of view -- by consideration of the $n$-dimensional vector space $\mathbb{F}_2^n$ over the field $\mathbb{F}_2=\{0,1\}$ which is an algebraic discrete structure. So, it can be a part of the Algebra module. We are convinced that its right place is in the DS module and most of the reasons justifying this opinion are indicated here and in the sample lecture. However, the lecturers in Algebra can extend the knowledge from the lecture by discussing the concepts of linearly dependence/independence of binary vectors, basis vectors, etc.

	The \textbf{primary audience} of this paper and the supporting lecture are the faculties of Computer Science (both lecturers and students) where the DS/DM module is studied. The secondary audience is the lecturers that teach subjects as Algebra, Programming in C (Java, C\#, etc.), Algorithms and Data structures, etc., and their students, as well as instructors looking for a student independent project, students in Master's degree programs, Ph.D. students, programmers, etc.
		
	This paper is organized as follows. In Section \ref{CS} we represent our case study about the $n$-dimensional Boolean cube. In Section \ref{CS_Concl} we summarize our observations and draw some conclusions after the case study. Section \ref{Lecture} starts with some comments on our sample lecture and description of its main parts and content. Thereafter, the methodology of teaching or (self)studying the lecture is discussed. Section \ref{Discussion} presents some arguments about the \textbf{benefits and results} of its study. Section 6 contains some concluding notes. 
\section{Case Study}
\label{CS}
	We examined how the $n$-dimensional Boolean cube and its properties are represented in some of the most popular textbooks, lecture notes (published on the Internet) and books in Discrete mathematics, Combinatorics, Algorithms and Data Structures, etc., most of which have several editions. We \textbf{explicitly note} that we do not try to evaluate these sources, most of which are excellent ones in their subject. We are interested in the selected topic only. Our study covers more than 80 such sources written in English, Russian and Bulgarian. A lot of them contain nothing about the Boolean cube and so they are not included in our list of references. Conversely, the famous book of Donald Knuth \cite{KNUTH} contains almost everything related to the $n$-dimensional Boolean cube. The book does not contain a special topic about it, but many terms, properties and applications of the cube are considered separately, in relation to many other discrete structures, data structures and algorithms. This is a lot more than we can represent here. The same is valid for the handbooks of Discrete mathematics (like the well-known handbook \cite{RMGGS}), as well as for the numerous papers that explore the Boolean cube from many different viewpoints. That is why such sources are not included in our case study. In the rest of the sources studied by us, the $n$-dimensional Boolean cube and related to it notions are considered in many different ways, varying widely from ``almost nothing" to ``almost complete". We can divide these textbooks and books into three main groups depending on the completeness of representation of these concepts, their structural properties, applications, etc. 
\subsection{Low Level}
\label{CS_LL}
	Many textbooks in this group just mention or use some notions as: \textit{bit strings, words of bits, binary sequences, binary words, binary representation} of integer, etc. \cite{AUDI, BP, DALE, GKP, GMTR, NR}. In addition, the concept of \textit{Gray code} (mostly in relation to \textit{Hamiltonian cycles} in graphs) is considered in \cite{JAND, GROSS, BSD, SDB}. Furthermore, \textit{Huffman code}, some other types of \textit{binary codes} and concepts as a \textit{weight} of a codeword, \textit{Hamming distance}  between codewords, \textit{minimal distance} of a code, etc., are discussed in \cite{AKIM, ANDI, CHEN, DPD, LPV, MERR, YABL}, etc. The \textit{bijection} between a subset of a given set and a binary number, represented in the necessary number of bits, is just mentioned in \cite{ERUS, LPV}. Haggarty \cite{HAGG} gives two examples where he defines informally the set of all sequences of zeros and ones (i.e., the bit strings) of length $n$. He uses the notion of \textit{characteristic vector} to show (informally again) the performance of basic operations on sets by logical operations on its characteristic vectors. However, the bijection between a subset and its characteristic vector is not mentioned. Almost the same is valid for the textbook \cite{IVBM}. Topics as \textit{Boolean algebra}, \textit{Boolean functions}, \textit{Boolean matrix algebra}, \textit{Boolean matrix multiplication}, etc. are considered in \cite{AKIM, DPD, ERUS, YABL}.
\subsection{Middle Level}
\label{CS_ML}
	The sources in this group consider the concepts related to the Boolean cube more comprehensively in comparison with these in the first group. Most of the textbooks (books) in this group are in the area of Algorithms and Data Structures. They accentuate to the usage of these concepts in their subject and in the practice. In addition to the binary representation of integers, topics as \textit{binary search, binary trees} (various types), \textit{binary heaps, files and dumps}, \textit{bit strings} (in relation to regular expressions), \textit{bitstreams}, as well as \textit{Huffman codes} and \textit{trees} (in relation to compression), \textit{Gray codes} (in relation to Hamiltonian paths and cycles), \textit{bitmaps}, \textit{bitwise operations}, etc., are considered in \cite{PNPD, SEDG1, SEDG2, SW4} and some other textbooks. The textbook \cite{PSSS} describes how to use certain functions from the Wolfram Mathematica in \textit{generating the subsets} of a given $n$-element set by binary strings of length $n$ -- in \textit{lexicographic order} (with respect to the binary strings or to the subsets), in a \textit{Gray code}, and also generating a random $k$-element subset. The authors comment on the efficiency of generating and show the relation between codewords in lexicographic order and these in Gray code. Functions for \textit{ranking} and \textit{unranking} of combinatorial objects are also discussed. The concept of an \textit{$n$-dimensional Boolean} cube is defined recursively and it is considered as a graph. Some of its combinatorial properties are given. The same topics, but quite more exhaustively are considered in \cite{RUSKEY}. The textbook of \cite{GGSA} starts with a section devoted to the \textit{$n$-dimensional Boolean cube}. Many important notions concerning it are defined. Some relations between them and the operations on binary vectors are explained and illustrated by examples. Furthermore, the Boolean cube is discussed in relation to \textit{Boolean functions}. Many solved problems and exercises are also given. 
\subsection{High Level}
\label{CS_HL} 
	The textbooks (books) in this group include more terms and properties besides those mentioned in the previous groups. Most of them are discussed in a lot more detail. Many authors consider the operations on vectors of the Boolean cube and relate them to operations on subsets of a given set. Some of them emphasize the structural properties. For example:
\begin{itemize}
\hyphenation{Boo-le-an}
\item Epp \cite{EPP} discusses in detail \textit{binary logic and arithmetic} -- digital logic circuits, Boolean expressions, binary representation of integers, one's and two's complement, the basic arithmetic operations on binary numbers, and also converting integers from/into decimal, binary, octal and hexadecimal numbers\footnote{For details, you can see the classic -- Chapter 4 of \cite{KNUTH2}.}. Almost the same topics are considered in \cite{GRIM, KOSHY, CLR, GRTJ}, etc.
	\item \textit{Computer representations of the subsets} of a given universal set $U$ and the \textit{basic operations on the subsets} are considered in relation to the \textit{bitwise operations} ``Not", ``And", ``Or", ``XOR" on  binary vectors of length equal to $|U|$. For this purpose, the notion of \textit{characteristic vector} (or \textit{characteristic function}) is defined and used in \cite{AUH1, KOSHY, KM}, or it is used informally in \cite{AUH2, GRIM, KUZNO, RND, ROS, CLR, SSKI}, etc. 
	\item \textit{Structural properties} -- the \textit{bijection} between the subsets of $U$ and the binary vectors of the $|U|$-dimensional Boolean cube is proven or shown informally in most of the textbooks, mentioned in this group. In \cite{GRTJ, EPP}, etc., \textit{algebraic structures} as Boolean algebras on sets are considered. Furthermore, the \textit{isomorphism} between these structures, or more generally, between finite Boolean algebras, are derived in \cite{GRIM, KUZNO, KM}, etc. 
	\item \textit{Relations} -- in general, and in particular as \textit{lexicographic precedence} and \textit{precedence} of binary vectors -- are introduced and their properties are discussed. Important notions as \textit{total} and \textit{partial orders} corresponding to these relations, and \textit{transitive closure} of relations\footnote{Many sources introduce a \textit{Boolean matrix multiplication} and apply it for computing  the transitive closure of relations or graphs.} are also discussed. The textbooks (books) in this group, as well as most of those in the previous groups, emphasize on the orders of the vectors of the Boolean cube \textit{in a Gray code} and their applications.  \textit{Generating the subsets} of a given set (or \textit{generating combinations}) -- in a lexicographic order, in a Gray code, etc. -- are considered in detail in \cite{RUSKEY, SSKI, RND, PSSS, GRIM, KOSHY, ROS}.  
	\item Many related topics, discrete structures and notions are also included -- for example, different types of \textit{binary codes} (linear, error-detecting, error-correcting, optimal prefix (Huffman) codes and trees, etc.), \textit{designs}, \textit{Boolean algebras}, \textit{Boolean functions}, \textit{Boolean matrices}, etc. Furthermore, the \textit{operations} on them and their structural properties are discussed. Topics from the area of \textit{hashing} (hash functions), \textit{cryptography} and \textit{data compression} are also considered in a lot of the mentioned textbooks. Different types of \textit{binary data structures} as binary trees, binary heaps, binary search trees, etc., and algorithms processing them are discussed in almost all examined textbooks in Algorithms and Data Structures. But not only there -- for example, Binary search, Merge sort, binary trees and their usage in representation of arithmetic expressions are explained and illustrated very well in \cite{EPP, ROS}.	
\end{itemize}
\section{Conclusions after the Case Study}
\label{CS_Concl}
	As we have noticed, the \textbf{main conclusion} after our case study is that the concepts related to the Boolean cube are important, but they are not included in a separated topic in all those textbooks. Having in mind only this topic, there are excellent sources (for example, \cite{GRIM, KOSHY, ROS}, etc.) that contain most of the notions concerning the Boolean cube and their properties, in comparison with those outlined in Section \ref{Lecture}. However, the general term that unifies them, i.e., the $n$-dimensional Boolean cube, is lacking. Hence a \textbf{separate topic about the Boolean cube is probably not included in the DS/DM module}. Some other observations and conclusions about the representation of the $n$-dimensional Boolean cube, the related notions, their properties and applications are:
\begin{enumerate}
	\item They vary widely depending on the authors' viewpoints, the subject of their textbook, its range, its size (the topic's content is represented more completely in the larger books), etc. 
	\item The level of formality in representation also varies widely -- from informal explanations or ideas for the notions (often given in examples, problems, etc.) to formal definitions, theorems, proofs, corollaries, etc.
	\item In the literature (including here and in the lecture) there are terms and notations that are not commonly accepted. Often these terms are considered separately and incompletely. Depending on the context, different terms and notations for the same notions are used even in the same textbooks. The connections between them, their properties and applications are given in part or omitted. That is why, it is difficult for the students to recognize them, as well as to understand that they all are simply binary vectors of certain lengths. What is more, they find it difficult to make any conclusions about them. Probably some techniques and tricks remain not enough understood or known for them -- for example, these for working with bit-wise representations of integers, the bit-wise operations on them, the properties of the different types of the vectors' orders etc.
	\item The case study confirmed our opinion that the discussed topic is important for many subjects, not only DS/DM, Programming, Algorithms and Data Structures. In addition to the related disciplines mentioned in the article, this topic is important for the training in Competitive Programming, as well as for the practice of future programmers.	
	\item The case study showed us what the most important notions, their properties, applications and relations to other subjects are. The last version of the lecture is based on those and some other observations. 
\end{enumerate}
\section{About the Lecture on the $n$-dimensional Boolean Cube and the Methodology of its Study}
\label{Lecture}
	The sample lecture about the Boolean cube is not included in the paper because of the following reasons:
\begin{itemize}
	\item The size of the lecture -- if we include it, the paper will be too large.
	\item It and its supporting teaching resources are available for free use (under the CC BY-SA 4.0 license) at \cite{VB_Lect}. We propose the .pdf and .tex files of the lecture, as well as the .pdf files of the figures. Thus the lecturer can exclude or modify some of its parts, or append additional parts in accordance with his/her viewpoint, the necessary time for its teaching, educational objectives, etc. The lecture contains many comments, remarks, references, etc., intended mainly to the lecturer -- they can also be excluded. Furthermore, the lecturer can easily prepare a presentation of the lecture.
	\item Possible violation of the publisher's copyrights.
\end{itemize}

		The concept of $n$-dimensional Boolean cube is introduced by the author as an example in the lecture ``Cartesian product of sets". It is defined as $\{0,1\}^n=\{(x_1, x_2,\dots, x_n)|\, x_i\in\{0,1\}, i=1, 2,\dots, n\}$ -- the set of all binary vectors. Only a few simple problems (for example, for counting binary vectors with fixed coordinates) can be solved by it. However, this is definitely not enough for our educational objectives and this first idea is upgraded by the lecture. It consists of 7 parts: 1.~Basic Notions, 2.~Relations and Orders (having 3 subsections), 3.~Operations on the Vectors of the Boolean Cube, 4.~Structural Properties of the Boolean Cube, 5.~Applications (having 3 subsections), 6.~Exercises (38), and 7.~Answers, Hints and Solutions. It contains more than 14 definitions, 8 theorems, 12 pictures, a lot of examples, remarks and comments, etc. Generally said, the lecture's content corresponds to the observations and conclusions we made in the previous section -- it includes the most important concepts, their properties, applications, etc. This dry and poor description of the lecture probably creates a wrong impression of it. Therefore, we strongly recommend the reader to get and read the original lecture.
	
	\textbf{The methodology of teaching} the lecture would depend on the amount of time the lecturer can devote to it in the lectures and labs.
The author uses (and proposes) two levels:
\begin{itemize}
	\item \textbf{Basic level}: includes teaching of sections 1. ``Basic notions", 2. ``Relations and orders", as well as a part of 5. ``Applications" (basically 5.1 ``Binary representation of integers, characters and strings") of the lecture. They can be taught in one lecture hour (of 45 minutes) by using a presentation. At least one lab hour is necessary for exercises\footnote{Almost all exercises after the lecture are appropriate for this goal.}, examples, details, etc. In this way, the students get basic knowledge on the topic. They obtain the ability to solve more enumeration problems, especially by modeling of coin tossing, encoding 2-directional paths (in particularly, at binary trees), sequences of switches with ON/OFF states, etc. The obtained knowledge is used in the remaining part of the DS/DM module -- for example, in solving enumeration problems in graphs, in illustrations (examples) of different types of graphs (Hamiltonian, Eulerian, bipartite, etc., as it is shown at the end of Section 2.2 of the lecture), in studying Boolean functions -- in defining and counting Boolean functions (especially those from the Post's classes), and so on. The basic knowledge is applicable in almost all related modules mentioned in the beginning of the paper.
\item \textbf{Extended level}: includes teaching the remaining parts of the lecture, which can be taught by a one-hour presentation. The last 4 exercises after the lecture can be done in a half lab hour or given for homework -- for developing programs that implement vector operations, generating the subsets of a given set, the basic operations on sets by vector operations (see \cite{AUH1, AUH2, SSKI, YORDZ}). This level extends and upgrades the previous knowledge by: operations on binary vectors, structural properties of the Boolean cube, the concept of characteristic vector, an example of an isomorphism between two discrete structures and between two partially ordered sets, more applications, etc. The new knowledge is related to the representation of graphs (the $i$-th row/column of the adjacency matrix of a directed graph is a characteristic vector of the set of successors/predecessors of the $i$-th vertex $v_i$), operations on graphs, etc. The bitwise operations on binary words are of great importance to the programmers, as it is shown in the next section. It is applicable to subjects such as Programming (in C, C\#, Java, etc.), Algorithms and Data structures, Cryptography, Competitive programming, etc. So, the knowledge obtained at this level gives a much deeper insight into the topic.
\end{itemize}

	When there is not enough time in the DS/DM module to teach the lecture or to do the exercises, the teacher may choose another studying method -- an \textbf{assignment of a self-study topic}. So, the students have to read the lecture or selected parts of it in order to solve some selected problems from the lecture, given as homework, or to work on a course project, related to the topic. In addition, students, lecturers, and practitioners who have a special interest in this topic can find some necessary knowledge in one lecture, instead of searching and reading a lot of books/textbooks. The lecturer should be ready to help the students in learning and understanding the lecture. In addition, if he/she wants to organize and lead the self-studying of the lecture, he/she can transform the studying into a tutorial format.
\section{Discussion: Is it Worth Studying the Topic about the Boolean Cube?}
\label{Discussion}
	The author is clear with some of the most popular objections about  including a new topic in the DS/DM (or some other) course. For example:
\begin{itemize}
	\item No evidence -- I expected to see the usual pedagogical experiment and analysis of its results to help me in my decision.
	\item There is a schedule for all topics in the course, it cannot be violated. Or, there is not time for lectures on new topics. Or, there are more important topics that remain outside the course.
	\item The politic and the curricula in our university are oriented to practice and to the requirements of the labor market.  The students need more and more practical knowledge and skills, they do not need more theoretical (unnecessary) knowledge. Some students even dislike mathematical subjects.
	\item The lecturers are free to choose what to teach as they understand their subject.
	\item The lecture does not correspond to my course -- it is too early for it, or it is too hard for my students. Maybe it is more appropriate for some master's degree programs in Computer Science. 
	\item I read the sample lecture, but I dislike it, etc.
\end{itemize}
	
	To answer the first objection we can say that such an experiment is unnecessary -- we teach new knowledge, with what should we compare it? So the results are known in advance. Obviously, they will confirm the benefits of studying the topic -- we claim this in view of everything presented here and our experience. Maybe some of the other objections are reasonable, but here we cannot argue with the reader who thinks in this way. We also believe in the known wisdom that if someone does not want to do something, he/she looks for reasons to reject it, otherwise, he/she looks for a way to do it. We hope the remaining part of this section will help the readers who have any doubts or hesitations. 
	
	It is our conviction that studying such a topic is worth it, and we present some arguments in support of it. Firstly, the conclusions in Section \ref{CS_Concl} show some obvious reasons for studying the discussed topic in the DS/DM module. Secondly, we hope the sample lecture is convinced enough. Besides these, we can point out other strong arguments about the benefits of its studying:
\begin{enumerate}
\item In \cite[p. 27]{CSC2013} it is noted: ``We expect curricula will have modules that incorporate topics from multiple knowledge areas", and on page 76: ``The material in discrete structures is pervasive in the areas of data structures and algorithms but appears elsewhere in computer science as well.", etc. The topic about the  Boolean cube satisfies these requirements completely. It incorporates topics from multiple knowledge areas and satisfies some other requirements from \cite{CSC2013}. To show this, it is enough to compare the knowledge areas from the case study and those areas of the subjects: Algorithms and Complexity, Architecture and Organization, Computational Science, Information Assurance and Security and especially Cryptography, Software Development Fundamentals, etc. Furthermore, the case study makes us believe that most of the lecturers in DS/DM are familiar enough with the knowledge areas of the subjects, related to DS/DM. We consider that they follow the recommendations in \cite{CSC2013} about the interdisciplinary connections. When it is possible, they illustrate or extend the knowledge obtained in previous subjects by appropriate relations, examples, explanations, etc., from their subject. They also give appropriate references for the forthcoming courses. Teaching a topic on the Boolean cube is another step in that direction.
\label{CSC}
	\item More reasons for studying the topic about Boolean cube can be found in the famous book \cite{KNUTH}. For example, in Section 7.1. Zeros and Ones, Knuth says: ``The amazing ability of 0s and 1s to encode information as well as to encode the logical relations between items, and even to encode algorithms for processing information, makes the study of binary digits especially rich.". Furthermore, in Section 7.1.3. Bitwise Tricks and Techniques, he notes: ``But we will see that Boolean operations on binary numbers deserve to be much better known. Indeed, they're an important component of every good programmer's toolkit.".
	\item Our practice and experience show that some notions, properties and applications concerning the Boolean cube that seem obvious to us, are not obvious to the students. That is why the studying of the considered topic in the DS/DM module is very useful. It becomes a base and we (or our colleagues) can refer to it always when we (or they) need. It includes good examples of already studied notions -- sets (subsets, power sets), relations and orders (POSets, chains, etc.), functions, combinatorial configurations; basic data types (as boolean, char, signed or unsigned short, int, long long); logical expressions; bitwise operations, etc. In this way, we recall, summarize and extend the obtained knowledge. Furthermore, the students study a lot of new notions, their properties and the connections between them; what they mean and how they can be used. 
	\item In \cite[p. 8]{GD} it is noted: ``Most exponential time algorithms are merely variations of exhaustive search, whereas polynomial time algorithms generally are made possible only through the gain of some deeper insight into the structure of a problem". The Boolean cube along with the operations on its vectors make it a discrete structure\footnote{It is very useful in the author's research, also for his colleagues working on similar research areas.}, hidden or standing behind many problems. Furthermore, the name of the module ``Discrete structures" requires the lecturer to show some discrete structures (there are not a lot of them in the module) and to drive the students' attention to their structural properties.	Their usage allows deriving more properties of the considered problems and accelerating the corresponding algorithms.
	\item	Benefits for the lecturers -- the new knowledge from the lecture can be transferred to teaching the courses Programming, Algorithms and Data structures, Competitive programming, etc., as it is shown in Section 5. ``Applications" of the lecture. For example, in the module on competitive programming, the author demonstrates the design and analysis of some algorithms in order to make them shorter and faster. He also discusses some problems from programming competitions where an analytical solution (i.e., a simple formula) can be derived -- it replaces time-consuming solutions, for example, of the type ``brute force'', ``backtracking'', etc. In addition, we have already commented on some connections of this knowledge with the Algebra course, as well as with Graph theory and Boolean functions. The reader of the lecture will understand why and how this knowledge can be transferred to subjects such as Combinatorial Algorithms, Cryptography, Data Compression, etc.
\item Benefits for the students -- besides the new notions, their properties, and applications, this knowledge:
\begin{itemize}
	\item is valuable, useful and important itself;
	\item proposes new tools and models of thinking in solving enumeration problems, in the design, analysis, and illustration of algorithms, in the understanding of the concepts -- already studied and forthcoming, etc.;
\hyphenation{Boo-lean}
	\item demonstrates the advantage of the usage of inductive and constructive definitions -- Definition 1.1 of the $n$-dimensional Boolean cube is used in the subsequent theorems in order to prove (by mathematical induction) some properties of the related concepts; 
	\item gives good examples of already studied notions (sets, operations on sets, relations, orders, partially ordered sets, etc.), discrete structures and isomorphism between them, etc.;
	\item can be applied to the forthcoming parts of the DS/DM course -- Graphs theory, Boolean functions, etc;
	\item connects, unites, ``sets bridges" between many subjects in the undergraduate programs in Computer Science;
	\item has many applications during the education and after it, in the practice of future programmers.
\end{itemize}
\end{enumerate}
\section{Conclusions}
\label{Concl}
	Here we represented the results of our examination on the completeness of representation of the $n$-dimensional Boolean cube in many books and textbooks. We shared our experience and we exposed our arguments about the necessity of studying such a topic in the DS/DM module. We consider that the notions, properties and applications concerning the Boolean cube should be taught together, in a lecture that relates and unites them. We hope that this paper will provoke the lecturers to think in this direction and make a positive decision. In such case, we also hope that the sample lecture and its supplementary materials (comments, references, problems and their answers, hints and solutions) will be useful to them. Finally, we hope that such a lecture will take its rightful place among the many free electronic resources in future textbooks.
\section*{Acknowledgement(s)}
%%
% \section*{Disclosure statement}
% No potential conflict of interest was reported by the author.
% \section*{Funding}
% \label{Funding}
 This work was partly supported by the Research Fund of the University of Veliko Turnovo (Bulgaria), Contract FSD-31-340-14/26.03.2019.

\begin{thebibliography}{99}

\bibitem{CSC2013} 
ACM/IEEE, Computer Science Curricula 2013: Curriculum Guidelines for Undergraduate Degree Programs in Computer Science. The Joint Task Force on Computing Curricula Association for Computing Machinery (ACM) and IEEE Computer Society (2013). Available at \\ http://ai.stanford.edu/users/sahami/CS2013/final-draft/CS2013-final-report.pdf. Last accessed 12 Feb. 2021

\bibitem{AUH1} 
Aho A.~V., Hopcroft J.~E., Ullman  J.~D., The Design and Analysis of Computer Algorithms. Addison-Wesley Publishing Company (1974).

\bibitem{AUH2} 
Aho A.~V., Hopcroft J.~E., Ullman J.~D., Data Structures and Algorithms. Addison-Wesley Publishing Company (1983).

\bibitem{ANDI} 
Anderson I., A first Course in Discrete mathematics. Springer-Verlag,  London (2001).

\bibitem{JAND} 
Anderson J., Discrete Mathematics with Combinatorics. Prentice-Hall, New Jersey (2001).

\bibitem{AKIM}
Akimov O.~E., Discrete mathematics. Logic, groups, graphs. Laboratoriya bazovih znanii, Moskow (2001). (in Russian)

\bibitem{AUDI}
Audibert P., Mathematics for Informatics and Computer Science. ISTE Ltd and John Wiley\&Sons, Inc. (2010).

\bibitem{VB_Lect}
Bakoev V., Shared lectures and files, ''An n-dimensional Boolean cube: basic notions, orders, properties and applications". The lecture is available at\\
\url{https://sites.google.com/view/valentin-bakoev/home}. Last accessed Dec. 20, 2021.

\bibitem{VB2014}
Bakoev V., Discrete mathematics: Sets, Relations, Combinatorics.  KLMN, Sofia  (2014). (in Bulgarian)

\bibitem{BSD}
Bogart K., Drysdale R., Stein Cl., Discrete Math for Computer Science Students (2004). Free online book, available at
http://www.freebookcentre.net/maths-books-download/Discrete-Math-for-Computer-Science-Students-(PDF-344P).html. Last accessed 12 Feb. 2021

\bibitem{BP} 
Brass P., Advanced Data Structures. Cambridge University Press. Cambridge (2008).

\bibitem{CHEN}
Chen W.~W.~L., Discrete Mathematics. (2008). Free online book, available at http://www.williamchen-mathematics.info/lndmfolder/lndm.html. Last accessed 12 Feb. 2021

\bibitem{CLR}
Cormen T., Leiserson Ch., Rivest R., Stein Cl., Introduction to Algorithms. 3rd edn. The MIT Press  (2009).

\bibitem{DALE}
Dale N., C++ Plus Data Structures. 3rd edn. Jones and Bartlett Publisher, Inc. (2003).

\bibitem{DPD}
Denev J., Pavlov R., Demetrovich Y.,  Discrete Mathematics. Nauka i izkustvo, Sofia (1984). (in Bulgarian)

\bibitem{EPP} 
Epp S., Discrete Mathematics with Applications. 4th edn. Brooks/Cole Cengage Learning (2010).

\bibitem{ERUS}
Erusalimski Y.~M., Discrete mathematics: theory, problems, applications. Vuzovskaia kniga, Moskow (2000). (in Russian)

\bibitem{GD}
Garey M., Johnson D., Computers and Intractability. A Guide to the Theory of NP-Completeness. W. H. Freeman and Company, New York (1979).

\bibitem{GRTJ} 
Garnier R., Taylor J., Discrete Mathematics for New Technology. 2nd edn. IOP Publishing Ltd. (2002).

\bibitem{GGSA} 
Gavrilov G., Sapozhenko A., Problems and Exercises in Discrete Mathematics, 3rd edn. FizMatLit, Moskow  (2005). (in Russian)

\bibitem{GMTR} 
Goodrich M., Tamassia R., Algorithm Design and Applications. John Wiley\&Sons Inc. (2015).

\bibitem{GKP}
Graham R., Knuth D., Patashnik O., Concrete Mathematics. A Foundation for Computer Science. 2nd edn. Addison-Wesley (1998).

\bibitem{GRIM} 
Grimaldi R., Discrete and Combinatorial Mathematics. An Applied Introduction. 5th edn. Addison-Wesley (2004).

\bibitem{GROSS} 
Gross J., Combinatorial Methods with Computer Applications. CRC Press (2008).

\bibitem{HAGG}
Haggarty R., Discrete Mathematics for Computing. Pearson Education Limited, UK (2001). 

\bibitem{HARARY}
Harary F., Hayes J.P., Wu H.-J., A Survey of the Theory of the Hypercube Graphs. Computers and Mathematics with Applications, 15 (4), 277-289 (1988).

\bibitem{IVBM}
Ivanov B.~M., Discrete mathematics. Algorithms and programs. Laboratoria Bazovih Znanii, Moskow (2003). (in Russian)

\bibitem{KNUTH2}
Knuth D., The art of computer programming, Volume 2:  Seminumerical algorithms. Addison-Wesley (1969).

\bibitem{KNUTH}
Knuth D., The art of computer programming, Volume 4A: Combinatorial Algorithms, Part 1. Addison-Wesley  (2011).

\bibitem{KOSHY} 
Koshy T., Discrete Mathematics with Applications. Academic Press (2003).

\bibitem{KUZNO}
Kuznetsov O., Discrete mathamatics for engineers. 6th edn.  Lan, St. Peterburg-Moskow-Krasnodar (2006). (in Russian)

\bibitem{LPV}
Lovasz L., Pelikan J., Vesztergombi K., Discrete Mathematics: Elementary and Beyond (Undergraduate Texts in Mathematics). Springer  (2003).

\bibitem{KM} 
Manev K., Introduction to Discrete Mathematics. 4th edn. KLMN, Sofia (2007). (in Bulgarian)

\bibitem{MERR}
Merris R., Combinatorics. 2nd edn. John Wiley\&Sons Inc. (2003).

\bibitem{PNPD}
Nakov P., Dobrikov P., Programming= ++ Algorithms. 3rd (revised) edn.  TopTeam Co, Sofia (2005). (in Bulgarian) Free online book, available at http://www.programirane.org/. Last accessed 12 Feb. 2021

\bibitem{NR}
Neapolitan R., Foundations of algorithms. 5th edn. Jones\&Bartlett Learning  (2015).

\bibitem{PSSS}
Pemmaraju S., Skiena S., Computational Discrete Mathematics: Combinatorics and Graph Theory with Mathematica. Cambridge University Press  (2003).

\bibitem{RND}
Reingold E., Nievergelt J., Deo N., Combinatorial algorithms, Theory and practice. New Jersey, Prentice-Hall  (1977).

\bibitem{ROS} 
Rosen K. H., Discrete Mathematics and its Applications. 7th edn. McGraw-Hill (2012).

\bibitem{RMGGS} 
Rosen K. (Editor in Chief), Michaels J., Gross J., Grossman J.,  Shier D., Handbook of Discrete and Combinatorial Mathematics. CRC Press (2000). 

\bibitem{RUSKEY}
Ruskey F., Combinatorial Generation. (Working Version (1j-CSC 425/520)  (2003). Available at http://page.math.tu-berlin.de/~felsner/SemWS17-18/Ruskey-Comb-Gen.pdf. Last accessed 12 Feb. 2021

\bibitem{CSavage}
Savage C.,  A Survey of Combinatorial Gray Codes. SIAM Review 39(4), 605--629 (1997).

\bibitem{SEDG1}
Sedgewick R., Algorithms. Addison-Wesley (1983).

\bibitem{SEDG2}
Sedgewick R., Algorithms in C. 3rd edn. Addison-Wesley (1998).

\bibitem{SW4} 
Sedgewick R., Wayne K., Algorithms. 4th edn. Addison-Wesley (2011).

\bibitem{SSKI} 
Skiena S., The Algorithm Design Manual. 2nd edn. Springer (2008).

\bibitem{SDB}
Stein Cl., Drysdale R., Bogart K.,  Discrete Mathematics for Computer Scientists. Addison-Wesley (2011).

\bibitem{YABL} 
Yablonski S., Introduction to discrete mathematics. 4th edn. Moskow: Visshaia Shkola (2003). (in Russian)

\bibitem{YORDZ}
Yordzhev K., The Bitwise Operations in Relation to the Concept
of Set, Asian Journal of Research in Computer Science, 1(4): pp. 1--8, 2018; Article no. AJRCOS.44314

\end{thebibliography}
\end{document}